\documentclass[12pt,a4paper,
]{article}
\usepackage{amsfonts,amssymb,latexsym}
\usepackage{theorem}

{\theorembodyfont{\rm}
  \newtheorem{defi}{Definition}[section]
  \newtheorem{rem}[defi]{Remark}
  \newtheorem{exa}[defi]{Example}
  \newtheorem{exas}[defi]{Examples}}
  \newtheorem{lem}[defi]{Lemma}
  \newtheorem{prop}[defi]{Proposition}
  \newtheorem{thm}[defi]{Theorem}
  
\newcommand{\Aa}{{\mathbb A}}

\newcommand{\PP}{{\mathbb P}}

\newcommand{\cB}{{\mathcal B}}
\newcommand{\cC}{{\mathcal C}}
\newcommand{\cD}{{\mathcal D}}

\newcommand{\fC}{{\mathfrak C}}

\newcommand{\End}{{\mathrm{End}}}

\newcommand{\dis}
                 {{\mathrel{\scriptstyle{\triangle}}}}

\newcommand{\GF}{{\mathrm{GF}}}
\newcommand{\GL}{{\mathrm{GL}}}
\newcommand{\PGL}{{\mathrm{PGL}}}

\newcommand{\M}{{\mathrm{M}}}
\let\phi=\varphi

\let\theta=\vartheta

\newcommand{\SDelimArray}[4]{\hbox{\scriptsize\arraycolsep=.5\arraycolsep
  $\left#1\!\!\begin{array}{*{#3}{r}}#4\end{array}\!\!\right#2$}}

\newcommand{\SMat}{\SDelimArray()}

\newcommand{\qed}{\ $\square$}
\newcommand{\pf}{{\sl{Proof:}\ }}
\date{\normalsize{\em Dedicated to
 H.\ Reiner Salzmann on the occasion of his 70th birthday.}}
\begin{document}
\title{Extending the Concept of Chain Geometry}
\author{Andrea Blunck\thanks{Supported by a Lise Meitner
 Research Fellowship
of the Austrian Science Fund (FWF), project M529-MAT.} \and
{Hans Havlicek}}

\maketitle

\begin{abstract} \noindent
We introduce the chain geometry $\Sigma(K,R)$
 over a ring $R$ with
a distinguished subfield $K$, thus extending the usual concept where $R$
has to be an algebra over~$K$.
A chain is uniquely determined by three of its points, if, and only if,
the multiplicative group of $K$ is normal in the group of units of $R$. This
condition is not equivalent to $R$ being a $K$-algebra. The chains through
a fixed point fall into {\em compatibility classes} which allow to describe
the residue at a point in terms of a family of affine spaces with a common set
of points.

\end{abstract}
\parskip .1cm
\parindent0cm
\section{Introduction}

Chain geometries $\Sigma(K,R)$, where $R$ is an
associative  algebra  over some  field~$K$,
have been investigated by many authors. They were first introduced
for the case of arbitrary commutative algebras by W.~Benz; compare his
monograph~\cite{benz-73}. Later, A.~Herzer and others considered also algebras that
are not commutative. For a survey, see \cite{herz-95}.

All these chain geometries are so-called {\em chain spaces}
(compare \cite{herz-95}). A chain
space is an incidence structure $\Sigma=(\PP,\fC)$, consisting of a
point set $\PP$ and a set $\fC$  of certain  subsets
of $\PP$, called {\em chains},   satisfying
the following axioms:

\begin{description}
\item[CS1] Each point lies on a chain,  and each chain contains at least
three points.
\item[CS2] Any three pairwise distant points lie together on exactly one chain.

Here two points $p,q\in\PP$ are called {\em distant}
(denoted by $p\dis q$), if they are different and joined
by at least one chain.
\item[CS3] For each point $p\in\PP$, the {\em residue} $\Sigma_p:=(\PP_p,
\fC_p)$, where
$\PP_p:=\{q\in\PP\mid q\dis p\}$ and
$\fC_p:=\{\cC\setminus\{p\}\mid p\in\cC\in\fC\}$,
 is a {\em partial affine space}, i.e., an incidence structure
resulting from an affine space by removing
certain parallel classes
of lines.
\end{description}

In this paper  we want to discuss how things alter if  the subfield
$K$ of the ring $R$ is   not necessarily contained in the center
of $R$, and hence  $R$ is not necessarily
a $K$-algebra. The chain geometries arising then will turn out to
fulfil certain weakened versions of the axioms above: In
general there is more than one chain joining three pairwise distant
points. The block set of a residue  $\Sigma_p$
is partitioned  into equivalence
classes such that for each class the whole  point set of $\Sigma_p$
together with the blocks  of this class form a partial affine space.

In the literature, so far only special cases have been discussed. See
 \cite{benz-73},IV\S2,  \cite{havl-93b}, \cite{havl-94a},
\cite{havl-97}, and \cite{maeu+m+n-80},
where  $K$ and $R$  are  skew fields, and
\cite{blunck-99b}, where $R$ is the ring of endomorphisms of a left
vector space over the skew field  $K$.

\section{Definition and basic results}

Let $K$ be a (not necessarily commutative) field, and let $R$ be a ring
with~$1$ such that $K\subseteq R$ and $1_K=1_R$. By
$R^*$ we denote the set of invertible elements (the {\em units})
 of $R$.
The ring $R$ is a $K$-algebra exactly if the field $K$ belongs to the
{\em center} $Z(R):=\{z\in R\mid \forall r\in R: rz=zr \}$ of $R$.

We  define the {\em chain geometry}
$\Sigma(K,R)$ over $(K,R)$  as this is done for algebras in
\cite{benz-73} and \cite{herz-95}.

The most important ingredient of the definition of $\Sigma(K,R)$
is the group
$$\Gamma:=\GL_2(R)$$
of invertible $2\times2$-matrices
with entries in $R$.
It acts in a natural way (from the right) on the free
left $R$-module $R^2$.

The {\em point set} of $\Sigma(K,R)$ is the {\em projective line}
over the ring  $R$. This is
the orbit $\PP(R):=R(1,0)^{\Gamma}$ of the free cyclic
submodule $R(1,0)\le R^2$ under the action of  $\Gamma$.
So
$$\PP(R)=\{R(a,b)\mid\exists c,d\in R:\SMat2{a&b\\c&d}\in \Gamma\}.$$
In other words, the elements of $\PP(R)$ are exactly those free cyclic
submodules of $R^2$ that possess a free cyclic complement.
Compare  \cite{herz-95} for  basic properties of
$\PP(R)$.

Note that in the special case that the ring $R$ is a field, this definition
coincides with the usual one, namely, $\PP(R)$ then is the set of all
$1$-dimensional subspaces of the vector space $R^2$.

Since we assume that $K$ is a subfield of $R$, the projective line
$\PP(K)$ over~$K$ can be embedded into $\PP(R)$  via
$K(k,l)\mapsto R(k,l)$. We call $\PP(K)$ (considered as a subset of
$\PP(R)$) the {\em standard chain} of $\Sigma(K,R)$, and denote it
by $\cC$.

The {\em chain set} of $\Sigma(K,R)$ is the orbit
$$\fC(K,R):=\cC^{\Gamma}.$$
 Altogether, the {\em chain geometry}
over $(K,R)$ is the incidence structure
$$\Sigma(K,R)=(\PP(R),\fC(K,R)).$$

By construction,  $\Sigma(K,R)$ satisfies
axiom {\bf CS1}  of  a chain space.

The kernel of the action of  $\Gamma=\GL_2(R)$  on $\PP(R)$
 is the center
$Z(\Gamma)$ of $\Gamma$, which coincides with
$\{\SMat2{z&0\\0&z}\mid z\in Z(R)^*\}$. So the group $\PGL_2(R)$
of permutations of $\PP(R)$ induced by $\Gamma$  is
isomorphic with $\Gamma/Z(\Gamma)$.
Since  $\fC(K,R)$ is an orbit under $\Gamma$,
the group $\PGL_2(R)$ consists of {\em automorphisms} of
$\Sigma(K,R)$.

Before investigating the incidence structure $\Sigma(K,R)$, we
first introduce the relation  `distant' on $\PP(R)$. It can be defined
on the projective line over any ring (cf.\ \cite{herz-95}).
Below we will see that in our
case it coincides with the relation `distant' defined  in  {\bf CS2}.
Note that some authors consider the relation `not distant' instead
and call it `parallel' (see, e.g., \cite{benz-73}).

Points $p=R(a,b)$ and $q=R(c,d)$ of $\PP(R)$ are called
{\em distant} ($p\dis q$), if the matrix
$\SMat2{a&b\\c&d}$ belongs to $\Gamma=\GL_2(R)$, i.e., if
$(a,b), (c,d)$ is a basis of $R^2$. Note that by this the relation
$\dis$ is well defined.

Just as  $\PP(R)$ and $\fC(K,R)$, also $\dis$ can be considered as an
orbit under $\Gamma$, namely,
\begin{equation}\label{Dis}\dis=(R(1,0),R(0,1))^\Gamma.
\end{equation}

Obviously, the relation $\dis$ is anti-reflexive and symmetric. Moreover,
if $R$ is a field, then $\dis$ equals the relation $\ne$.

This leads us to a characterization of $\dis$ in terms of the chain
geometry $\Sigma(K,R)$ (see \cite{herz-95}, 2.4.2, for the case of
algebras):

\begin{lem}\label{distant}
Let $p,q\in\PP(R)$ be different points of $\Sigma(K,R)$.
Then $p\dis q$ holds exactly if there is a chain $\cD\in\fC(K,R)$
joining $p$ and $q$.
\end{lem}
\pf
By (\ref{Dis}) we know that
 $p\dis q$ implies  $p=R(1,0)^\gamma$, $q=R(0,1)^\gamma$ for some
$\gamma\in \Gamma$. Hence in this case
$p,q\in\cC^\gamma\in\fC(K,R)$.

Conversely, if $p,q\in\cC^\gamma\in\fC(K,R)$ (with $\gamma\in\Gamma$),
then $p^{\gamma^{-1}}$
and $q^{\gamma^{-1}}$ are different points of $\cC=\PP(K)$.
 On  $\cC$ one has two `distance' relations:
The ordinary one on $\PP(K)$ (which is the relation $\ne$), and the
one inherited from $\PP(R)$. However, one can easily check that
the two relations coincide, because
$\GL_2(K)=\M(2\times2,K)\cap\GL_2(R)$,
where  $\M(2\times2,K)$ denotes the
ring of $2\times 2$-matrices over $K$. So we have $p^{\gamma^{-1}}\dis
q^{\gamma^{-1}}$. Since $\gamma$ preserves $\dis$, this proves the
assertion.
\qed

 We now want to determine the
chains through three given pairwise distant points. Note that, by
Lemma \ref{distant}, any two points on a chain must be distant.

\begin{prop}\label{ThreePoints}
Let $p,q,r\in \PP(R)$ be pairwise distant. Then there is at least one
chain $\cD\in\fC(K,R)$ containing $p$, $q$, and $r$.
\end{prop}
\pf The group $\Gamma$  acts
{\em $3$-$\dis$-transitively}   on $\PP(R)$, i.e., transitively on the set
of triples of pairwise distant points of $\PP(R)$ (see \cite{herz-95},
1.3.3). So there exists a $\gamma\in\Gamma$ with $p=R(1,0)^\gamma$,
$q=R(0,1)^\gamma$, $r=R(1,1)^\gamma$, and $\cD:=\cC^\gamma$ is a chain through
$p$, $q$, and $r$.
\qed

This means that $\Sigma(K,R)$ fulfils the existence part of
axiom {\bf CS2}. The uniqueness statement of   {\bf CS2},
however, will not hold in general.

The essential result on the group action of $\Gamma$ on $\Sigma(K,R)$
is as follows:

\begin{thm}\label{Neu}
Let $\cD,\cD'\in\fC(K,R)$, and let $p,q,r\in\cD$
and $p',q',r'\in\cD'$ be three pairwise distant points, respectively.
Then there exists a $\gamma\in\Gamma$ such that
$p^\gamma=p'$, $q^\gamma=q'$, $r^\gamma=r'$, and  $\cD^\gamma=\cD'$.
\end{thm}
\pf
There exists a $\gamma_1\in\Gamma$ mapping $\cD$ to the standard chain
$\cC$. Put $p_1:=p^{\gamma_1}$, $q_1:=q^{\gamma_1}$,
$r_1:=r^{\gamma_1}$. The group $\GL_2(K)\le\Gamma$ leaves $\cC$
invariant and acts triply transitively on $\cC$. Hence there is a
$\gamma_2\in\GL_2(K)$ with $p_1^{\gamma_2}=R(1,0)$,
$q_1^{\gamma_2}=R(0,1)$, $r_1^{\gamma_2}=R(1,1)$ (and $\cC^{\gamma_2}=\cC$).
Define $\gamma_1'$ and $\gamma_2'$ accordingly.
Then $\gamma=\gamma_1\gamma_2\gamma_2'^{-1}\gamma_1'^{-1}$ has
the required properties.
\qed

\begin{thm}\label{Thm1}
Let $\Sigma=\Sigma(K,R)$ and let
\begin{equation}\label{DefN} N:=N_{R^*}(K^*)=\{n\in R^* \mid n^{-1}Kn=K\}
\end{equation}
be the
 normalizer of $K^*$ in $R^*$.

Then the  set of chains
through a triple of pairwise distant points of $\Sigma$
 is in $1$-$1$-correspondence with the set
 $$R^*/N:=\{Nr\mid r\in R^*\}$$
of right cosets of $N$.

In particular, in $\Sigma$ there exists exactly one chain through
each triple of pairwise distant points if,
and only if, $K^*$ is normal in $R^*$.
\end{thm}
\pf The subgroup
\begin{equation}\label{Omega} \Omega:=\{\SMat2{a&0\\0&a}\mid a\in R^*\}\cong R^*
\end{equation} of $\Gamma$ is the stabilizer of the triple
$(R(1,0), R(0,1), R(1,1))$ of  {\em standard points}.
So, by Theorem \ref{Neu} the chains through
 $R(1,0)$, $R(0,1)$, $R(1,1)$
are exactly
the images $\cC^\omega$, $\omega\in \Omega$.
Since the stabilizer of $\cC$ in $\Omega$ is
\begin{equation}\label{OmegaC}
\Omega_\cC=\{\SMat2{n&0\\0&n}\mid n\in N\}\cong N,
\end{equation}
the assertion follows for the standard points and
thus, by Theorem \ref{Neu},
for any three pairwise distant points.
 \qed

So axiom  {\bf CS2} holds in $\Sigma(K,R)$ exactly if
$R^*=N$.
Of course, the condition $R^*=N$ is satisfied if $K$ belongs to the
center of $R$, i.e., if $R$ is a $K$-algebra. Hence Theorem \ref{Thm1}
reconfirms that  {\bf CS2} is valid for
chain geometries over algebras (compare Section 1).

However,
$R^*=N$ does not necessarily mean that $K$ is central in $R$.
We give two examples:

\begin{exas}\label{Exa1}
 \begin{enumerate}
 \item Let $K$ be a non-commutative field, and let $R=K[X]$ be
 the polynomial ring over $K$ in the central indeterminate $X$.
 Then $R^*=K^*$, and hence $N=R^*$. However, $K\not\subseteq
 Z(R)=Z(K)[X]$. \\
 By \ref{Thm1}, in  $\Sigma(K,R)$ there is  exactly one chain through
 any three given pairwise distant points.
 But here the situation is even more special because of the lack of
 units outside $K^*$: Any {\em two} distant points are joined by {\em exactly
 one} chain since the stabilizer of the pair $(R(1,0),R(0,1))$
 belongs to $\GL_2(K)$.
 \item Let $K=\GF(4)$, the field with $4$ elements. Then
 $K$ can be represented by certain $2\times2$-matrices
 over $F=\GF(2)$, i.e., we may assume $K\subseteq R=\M(2\times2,F)$
 (with $1_K=1_R$). Of course, $Z(R)=F\not\supseteq K$.
 Since  $|R^*|=6$, the group $K^*$ has index $2$ in $R^*$ and
 hence is normal. This means that  $N=R^*$.\\
  Also here the situation is rather  special: The chains are exactly the
  maximal sets of pairwise distant points.
 \end{enumerate}
\end{exas}

\begin{rem}\label{CBH}
In the following special
  cases  the validity of $R^*=N$ implies centrality
of $K$:
 \begin{enumerate}
 \item Let $K$ and $R$ be (not necessarily commutative) fields
 such that $R^*=N$ and $K\ne R$. Then $K\subseteq Z(R)$
 (Cartan-Brauer-Hua, see \cite{benz-73}, p.\ 323).
 \item Let $U$ be a left vector space  over $K$, and let $R=\End_K(U)$
 be the ring of endomorphims of $U$. We  embed $K$ into $R$ with
 respect to a fixed basis $(b_i)_{i\in I}$ of $U$ via
 $k\mapsto\lambda_k:(b_i\mapsto kb_i)$.

 If now $R^*=N$ holds, and $R\ne K$ (i.e., $\dim U>1$),
  then $K\subseteq Z(R)$.
 This  can be shown by calculating $\phi^{-1}\lambda_k
 \phi$ for certain elementary transvections $\phi\in R^*$
 (see \cite{blunck-99b}, proof of 4.1).
 In particular, in this case
 $K$ is  commutative and so $K=Z(R)$.
 \end{enumerate}
\end{rem}

Since the kernel of the canonical epimorphism $\GL_2(R)\to\PGL_2(R)$
consists exactly of the elements of $\Omega$ (see (\ref{Omega}))
where
the scalar $a$ belongs to $Z(R)^*=Z(R)\cap R^*$,
we have  the following:

\begin{rem}
The action of $\PGL_2(R)$ on $\PP(R)$ is {\em sharply} $3$-$\dis$-transitive,
exactly if the multiplicative group $R^*$ is contained in the center
$Z(R)$.
\end{rem}

Note that  the condition $R^*\subseteq Z(R)$ does not imply
that $R$ is commutative. Consider, e.g., the polynomial ring 
$K[X,Y]$ over a commutative field $K$
in the {\em non-commuting} indeterminates $X$ and $Y$.

We  turn back to the example given in \ref{Exa1}(b).
It can be generalized to arbitrary quadratic extensions of not necessarily
commutative
fields (cf.\ \cite{cohn-95}, Section 3.6).
We restrict ourselves to the finite case. Then we
can count the chains through  three pairwise distant points:

\begin{exa}\label{Exa2}
Let $q$ be any prime power $\ne 1$, and let
$F=\GF(q)$ be the   field with $q$ elements. Moreover, let $K=\GF(q^2)$
and $R=\M(2\times 2,F)$.
Then $K=F+Fi$ ($i\in K\setminus F$) with $i^2=s+ti$ (for suitable
$s,t\in F$). The right regular representation of $K$ yields an
embedding of $K$ into $R$, namely, $a+bi\mapsto
\SMat2{a&b\\bs&a+bt}$.

We now consider the elements of $R$ also as endomorphisms of the
vector space $K^2\supseteq F^2$.
A matrix in $R$ describes an element of $K$, i.e., it
has the form $ \SMat2{a&b\\bs&a+bt}$,  exactly if $(-i,1)$ is
one of its eigenvectors (see \cite{havl-94a}, Theorem~1).
The second eigenvector then  must be $(-\bar i,1)$,
where $\bar i$ is the {\em conjugate} of $i$ w.r.t. the
Galois group of
$K/F$.
So $N$ is the subgroup of $R^*$ that leaves  the set
$\{K(-i,1),K(-\bar i,1)\}$ invariant.

Now the group $R^*$ is the product of $K^*$ with the
subgroup $\{\SMat2{x&0\\y&1}\mid x\in F^*, y\in F\}$, which acts sharply
transitively on the set of vectors $(z,1)$,  $z\in K\setminus F$.
We conclude that $N=K^*\cdot\langle \kappa\rangle$, where $\kappa$ is the
unique matrix of type $\SMat2{x&0\\y&1}$ mapping $(-i,1)$ to $(-\bar i,1)$.
Obviously, $\kappa$ is an involution. Hence $|N|=2|K^*|=2(q^2-1)$.
Since $|R^*|=(q^2-1)(q^2-q)$, this means that in $\Sigma(K,R)$ there are
exactly $${1\over2}(q^2-q)$$ chains through three pairwise distant points.
In case $q=2$ this number equals $1$, as asserted in \ref{Exa1}(b).
\end{exa}

We now determine the intersection of all chains through three
pairwise distant points of a chain geometry $\Sigma(K,R)$.

\begin{prop}
Let $p,q,r\in\PP(R)$ be pairwise distant.
Then the
intersection of all chains through $p,q,r$ is an {\em $F$-chain},
i.e., the image of the projective line $\PP(F)$  over the subfield
\begin{equation}\label{DefF}
F:=\bigcap\limits_{u\in R^*}u^{-1}Ku
\end{equation}
of $K$ under a suitable $\gamma\in\Gamma$.
\end{prop}
\pf
We consider w.l.o.g.\ the standard points
$R(1,0)$, $R(0,1)$, $R(1,1)$. The chains joining them are exactly
the images $\cC^\omega$,
$\omega\in\Omega$ (compare (\ref{Omega})).
We compute
$$\bigcap\limits_{\omega\in\Omega}\cC^\omega
=\{R(1,0)\}\cup\bigcap\limits_{u\in R^*}\{R(u^{-1}ku,1)\mid k\in K\},
$$
which  equals $\PP(F)$, considered as a subset of $\PP(R)$.
\qed

Of course, in (\ref{DefF})
it suffices to let $u$ run over a system of representatives
for $R^*/N$.  In particular, if $R^*=N$, then $F=K$. This is clear also
because $R^*=N$ means that the chain through $p,q,r$ is unique.
In case $F\ne K$ the theorem of Cartan-Brauer-Hua (see \ref{CBH}(a)),
applied to $K$, implies $F\subseteq Z(K)$.
Moreover, $F$ always contains the subfield $K\cap Z(R)$ of $K$.
In certain  cases, $F$ and $K\cap Z(R)$ coincide:

\begin{exas}
 \begin{enumerate}
 \item Let $R$ be a skew field and $R\ne K$. Then $F\subseteq Z(R)$ by the
 theorem of Cartan-Brauer-Hua  and hence $F=K\cap Z(R)$.
 \item Let $K$ and $R$ be as in \ref{Exa2}. Then for any $u\in R^*
 \setminus N$ the field  $K\cap u^{-1}Ku$ is a proper subfield
 of $K=\GF(q^2)$ and hence equals $\GF(q)=Z(R)$.
 \item Let $R=\End_KU$ for some left vector space $U$ over $K$
 with $\dim U>1$.
 Then $F=Z(K)=Z(R)$ (compare \cite{blunck-99b}).
 \end{enumerate}
\end{exas}

\section{Compatibility of chains}

Now we want to investigate the set of all  chains through a fixed point.
For $p\in\PP(R)$ let
\begin{equation}\label{PPp}
\PP_p:=\{q\in\PP(R)\mid q\dis p\}
\end{equation}
and
\begin{equation}
\fC(p):=\{\cD\in\fC(K,R)\mid p\in\cD\}.
\end{equation}
We are going to introduce an equivalence relation called {\em compatibility}
on $\fC(p)$.
Since the group
$\PGL_2(R)$ acts transitively on $\PP(R)$ and consists of automorphisms
of $\Sigma(K,R)$, we may restrict ourselves to the case $p=R(1,0)$. We shall
denote this point also by the symbol $\infty$.

With \ref{Neu}
we obtain
that the  set $\fC(\infty)$ consists exactly of the images of
the standard chain $\cC$  under the group
\begin{equation}
\Gamma_\infty=\{\SMat2{a&0\\c&d}\mid
a,d\in R^*, c\in R\},
\end{equation}
 which is the stabilizer of $\infty$ in $\Gamma=
\GL_2(R)$.

The chains $\cB,\cD\in\fC(\infty)$  are called  {\em compatible
 at $\infty$}
 (denoted by $\cB\sim_\infty \cD$), if
 they belong to the same orbit under the action of
 the group
\begin{equation}
\Delta:=\{\SMat2{a&0\\c&1}\mid a\in R^*, c\in R\}\unlhd\Gamma_\infty
\end{equation}
 on $\fC(\infty)$.

By definition, compatibility is an  equivalence relation on $\fC(\infty)$.
Since $\Delta$ is  normal in $\Gamma_\infty$, compatibility is
invariant under the action of $\Gamma_\infty$.

The  equivalence classes w.r.t. $\sim_\infty$ are called
{\em compatibility classes}, the compatiblity class of $\cB\in\fC(\infty)$
is denoted by $[\cB]_\infty$.

For an element $\delta=\SMat2{a&0\\c&1}\in\Delta$ and a point $p=R(x,1)
\in\PP_\infty$ we compute $p^\delta=R(xa+c,1)$. This shows that the action
of $\Delta $ on $\PP_\infty$ is {\em sharply $2$-$\dis$-transitive}.

The following theorem is essential:

\begin{thm}\label{ThmComp}
Let $\cD,\cD'\in\fC(\infty)$, and let $p,q\in \cD\setminus\{\infty\}$
and $p',q'\in\cD'\setminus\{\infty\}$
be different points, respectively. Moreover, let $\delta$
be the unique element of $\Delta$ with  $p^\delta=p'$ and $q^\delta=q'$.
Then $\cD\sim_\infty\cD'$ holds exactly if $\cD'=\cD^\delta$.
\end{thm}
\pf
Let $\cD\sim_\infty\cD'$. Then $\cD'=\cD^{\delta'}$ for some $\delta'\in
\Delta$. Since the group
 $\Delta(K):=\Delta\cap \GL_2(K)$
 acts $2$-transitively on $\cC\setminus\{\infty\}$, there is a
subgroup of $\Delta$ (conjugate to $\Delta(K)$) acting $2$-transitively
on $\cD\setminus\{\infty\}$. So we may w.l.o.g. assume $p^{\delta'}=p'$,
$q^{\delta'}=q'$. Uniqueness of $\delta$ implies $\delta'=\delta$
and hence $\cD'=\cD^\delta$.\\
The  proof of the converse is obvious.
\qed

\begin{thm}\label{CorComp}
Let $\cD\in \fC(\infty)$,  and
let $p',q'$ be arbitrary distant points of $\PP_\infty$.
Then there is a unique chain $\cD'\in\fC(\infty)$ with $\cD'\sim_\infty\cD$
and $p',q'\in\cD'$.

In particular, each compatibility class in $\fC(\infty)$  has a unique
 representative through the standard points.
\end{thm}
\pf
Choose different points $p,q\in\cD\setminus\{\infty\}$. There
is a unique $\delta\in\Delta$ such that $p^\delta=p'$,
$q^\delta=q'$. By Theorem \ref{ThmComp},
 $\cD':=\cD^\delta$ is the only chain with the required
properties.
\qed

So the set $\{\cC^\omega\mid\omega\in\Omega\}$ (cf.\ (\ref{Omega}))
of all chains through the standard points is a complete set of
representatives for the compatibility classes in $\fC(\infty)$.
By \ref{Thm1}, this means that  $\fC(\infty)/{\sim_\infty}=
 \{[\cB]_\infty\mid \cB\in\fC(\infty)\}$ is in $1$-$1$-correspondence
with  $R^*/N$.

We shall need the following explicit description of the chains
compatible at $\infty$
with the standard chain $\cC$:

\begin{lem}\label{LemComp}
Let $\cD=\cC^\gamma\in\fC(\infty)$, where
$\gamma=\SMat2{a&0\\c&d}\in\Gamma_\infty$.
Then  $\cD\sim_\infty\cC$ is equivalent to  $d\in N$.
\end{lem}
\pf
The unique
$\delta\in\Delta$ with $R(0,1)^\delta=R(0,1)^\gamma$ and
$R(1,1)^\delta=R(1,1)^\gamma$ equals
$\SMat2{d^{-1}a&0\\d^{-1}c&1}$. By Theorem \ref{ThmComp}, we have $\cC
\sim_\infty\cC^\gamma$ if, and only if, $\cC^\gamma=\cC^\delta$,
or, in other words, if $\gamma=\omega\delta$ with
$\omega=\SMat2{d&0\\0&d}\in\Omega_\cC$.
This in turn is equivalent to $d\in N$ (cf. (\ref{OmegaC})).
\qed

The relation of compatibility can be
carried over to the set $\fC(p)$ of chains through an arbitrary
point $p\in\PP(R)$  in a natural way: We say that $\cB,\cD\in\fC(p)$
are {\em compatible at $p$} exactly
if $\cB^{\gamma}\sim_\infty \cD^{\gamma}$ holds  for $\gamma\in \Gamma$ with
$p^\gamma=\infty$. This is independent of the choice of $\gamma$
because $\Delta$ is  normal
in $\Gamma_\infty$.

Theorem \ref{CorComp} implies that any two different
chains with at least three common points are
{\em non-compatible} at each point of intersection.

For chains meeting only in two points the situation is different.
We  study the  chains containing $\infty$ and $0:=R(0,1)$.

\begin{prop}
Let $\cD$ be a chain through $\infty$ and $0$.
Then $\cD\sim_\infty\cC$ and $\cD\sim_0\cC$ holds
exactly if $\cD=\cC^\delta$, where
$\delta=\SMat2{a&0\\0&1}$ for some $a\in N$.
\end{prop}
\pf
Let $\cD$ be a chain through  $\infty$ and $0$ compatible at $\infty$ with
$\cC$.
By Theorem \ref{ThmComp}, we may assume
$\cD=\cC^\delta$ for  $\delta=\SMat2{a&0\\0&1}$ (with $a\in R^*$).
The chains $\cC$ and $\cD$ are
compatible at $0$ exactly if
$\cC^{\gamma}\sim_\infty\cD^{\gamma}$ holds for
some $\gamma\in \Gamma$ mapping $0$ to $\infty$. We choose
$\gamma=\SMat2{0&1\\1&0}=\gamma^{-1}$, and hence obtain
$\cC^{\gamma}=\cC$ and
$\cD^{\gamma}=\cC^{\gamma\delta\gamma}$. Since
$\gamma\delta\gamma=\SMat2{1&0\\0&a}$, we conclude from  \ref{LemComp}
that $\cC\sim_0\cD$
holds exactly if $a\in N$.
\qed

So, whenever there is more than one compatibility class at $\infty$
(i.e., in case $R^*\ne N$), one can find chains through $\infty$ and
$0$ compatible at $\infty$ and not compatible at $0$. Of course
the same holds for every other pair of distant points.
In particular, compatibility cannot be considered as
a global equivalence relation on the whole chain set.

\section{The residue at a point}

We define the residue at a point of
the chain geometry  $\Sigma(K,R)$ exactly as (for arbitrary
incidence structures) in
Section\ 1.

For a point $p\in\PP(R)$ we consider the point set
 $\PP_p$ as defined in (\ref{PPp}) and the {\em block
set} $$\fC_p:=\{\cD\setminus\{p\}\mid\cD\in\fC(p)\}.$$
The incidence structure $$\Sigma_p:=(\PP_p,\fC_p)$$ is the
{\em residue of $\Sigma:=\Sigma(K,R)$ at $p$}.

Again we may restrict ourselves to the case $p=\infty$. Each residue of
$\Sigma$ is isomorphic to $\Sigma_\infty$.

We compute $\PP_\infty=\{R(x,1)\mid x\in R\}$. We often identify the
point $R(x,1)$ with the element $x\in R$, and thus the set $\PP_\infty$
with $R$.

Next we investigate the block set $\fC_\infty$. We introduce
the {\em standard block} of $\Sigma_\infty$, this is $C:=\cC
\setminus\{\infty\}$.

The relation of compatibility  is carried over to $\fC_\infty$ from the
set $\fC(\infty)$: Two blocks $B,D\in\fC_\infty$ are called
{\em compatible} (at $\infty$), if the chains
$B\cup\{\infty\}$ and $D\cup\{\infty\}$ are compatible at $\infty$.

The compatibility class of the block
$B\in\fC_\infty$ is written as~$[B]_\infty$.
We are going to study the incidence structure $(\PP_\infty,[B]_\infty)$.
Because of \ref{CorComp} we may always assume
$B=C^\omega$ with  $\omega\in\Omega$.

\begin{thm}\label{IsoAffin}
Let $u\in R^*$, $\omega=\SMat2{u&0\\0 &u }$, $B=C^\omega$,
and $\alpha:x\mapsto u^{-1}xu$. Then the following statements hold:
\begin{enumerate}
 \item
 The fields $K$ and $u^{-1}Ku$ are isomorphic subfields
 of $R$. The associated affine spaces $\Aa(K,R)$ and $\Aa(u^{-1}Ku,R)$
 are isomorphic via the semilinear bijection $\alpha$.
  \item The incidence structure
$(\PP_\infty,[B]_\infty)$
is a partial affine space in the affine space $\Aa(uKu^{-1},R)$.
 \item
  The  isomorphism $\alpha:\Aa(K,R)\to \Aa(uKu^{-1},R)$
 of affine spaces induces the isomorphism
 $\omega|_{\PP_\infty}:(\PP_\infty,[C]_\infty)\to(\PP_\infty,[B]_\infty)$
 of partial affine spaces.
 \end{enumerate}
\end{thm}
\pf
(a): This is a straightforward calculation.

(b):
The  compatibility class $[B]_\infty$
consists of all sets
$$\{R(kua+uc,u)\mid k\in K\}=(u^{-1}Ku)a+c,$$
where $a\in R^*, c\in R$.
Hence  the blocks of $ \fC_\infty$  compatible
with  $B$ are certain lines of the affine
space $\Aa(u^{-1}Ku,R)$. More exactly, a line $u^{-1}Kux+y$
($x,y\in R$, $x\ne 0$)
of $\Aa(u^{-1}Ku,R)$ is a block of $[B]_\infty$ if, and only if,
$x$ is a unit in  $R$.

(c) follows from (a) and (b).
\qed



Because of \ref{LemComp}, the compatibility class $[C]_\infty$ equals
 the entire block set $\fC_\infty$
of the residue $\Sigma_\infty$
exactly if  $R^*=N$.
So in this case the residue $\Sigma_\infty$ (and thus also every
other residue of $\Sigma$) is a partial affine space, i.e., axiom {\bf CS3}
holds in $\Sigma$.

Together with \ref{Thm1} we have

\begin{thm}\label{ChainSpace}
The chain geometry $\Sigma(K,R)$ is a chain space if, and only if,
the multiplicative group $K^*$ is normal in $R^*$.
\end{thm}

The examples of \ref{Exa1} yield chain spaces $\Sigma(K,R)$
where $R$ is not a $K$-algebra. For  the chain space of \ref{Exa1}(b)
one can even show even more:

\begin{prop}\label{NonClassical}
Let $K=\GF(4)$ and $R=\M(2\times 2,\GF(2))$.
Then  the chain space $\Sigma=\Sigma(K,R)$ is
not isomorphic to any chain geometry
$\Sigma(L,S)$ over some $L$-algebra $S$.
Moreover, $\Sigma$ cannot be embedded into any chain geometry
 $\Sigma(L,S)$ over a {\em strong} $L$-algebra $S$ as a
{\em subspace}.

{\em
(For the definition of a subspace cf.\ \cite{herz-92}
and \cite{kroll-92}, for the definition of a strong algebra,
 cf.\ \cite{herz-92}.)}
\end{prop}
\pf
Assume first that  $\Sigma\cong \Sigma(L,S)$ for some $L$-algebra $S$.
Then $|\cC|=5$ implies $L\cong K=\GF(4)$, and $|\PP_\infty|=|R|=16$ implies
$\dim_LS=2$. There are  three types of $2$-dimensional algebras over
$L$, and the residues of the
associated chain geometries are obtained by removing at most
two  parallel classes of lines from the affine plane
over $L$ (see \cite{benz-73} or \cite{herz-95}). However, in
$\Sigma_\infty$ there are only two blocks through $0$, each containing
$3$ of the $6$ elements of~$R^*$. This means that {\em three} parallel
classes of lines of the affine plane are missing, a contradiction.

Now assume that $\Sigma$ is isomorphic to a subspace $\Sigma'$ of some
$\Sigma(L,S)$, where $S$ is a strong $L$-algebra.
Again we have $L\cong K=\GF(4)$.
By \cite{herz-92}, Theorem 2, the subspace $\Sigma'$ can be described by a
$2$-dimensional vector subspace $J$ of  $S$ (considered as a vector space
over~$L$)
which  contains $1$ and is closed w.r.t. squaring. This
already implies
that $J$ itself is  an $L$-algebra, and $\Sigma\cong \Sigma(L,J)$
contradicts the first assertion.
\qed

In the second part of this proposition we had to restrict ourselves to
the case of strong algebras becauses otherwise the proof of the
coordinatization theorem for subspaces (\cite{herz-92}, Theorem 2)
does not work. Since there seems to be no algebraic description
of the subspaces of the chain geometry over an arbitrary $L$-algebra $S$,
it remains open  whether the chain space $\Sigma=\Sigma(K,R)$
can be embedded into some $\Sigma(L,S)$ or not.

By Theorem \ref{IsoAffin},  the set $\PP_\infty=R$ is the common point set
of the family  $\Aa(u^{-1}Ku,R)$ of  affine spaces, where $u\in R^*$.

If we fix one such affine space, say $\Aa(K,R)$, then some
blocks of $\fC_\infty$  are lines of that space.
In the following Proposition, we  describe all  blocks,
using the fact that each block is a line in some affine
space $\Aa(u^{-1}Ku,R)$.

\begin{prop}\label{uK}
Let $u\in R^*$  and let
$B$ be a block appearing as a line of $\Aa(u^{-1}Ku,R)$.
Then the trace space induced on
$B$ by $\Aa(K,R)$ is isomorphic to the affine space $\Aa(F_u,K)$,
where $F_u$ is the subfield $F_u:=K\cap uKu^{-1} $ of $K$.
\end{prop}
\pf
The translation groups of $\Aa(K,R)$ and $\Aa(u^{-1}Ku,R)$
are the same. Hence
we may  assume that $B$ contains $0$. So  $B=u^{-1}Kua$
for some  $a\in R^*$.
One easily checks that $B$ is a left vector space over $F_{u^{-1}}
=K\cap u^{-1}Ku$. Moreover, $K$ is a left vector space over
$F_u$, and $\alpha:x\mapsto u^{-1}xu$ is an isomorphism
$F_u\to F_{u^{-1}}$.
The mapping  $\iota:k\mapsto u^{-1}kua=k^\alpha a$ is a semilinear
bijection $K\to B$
with accompanying  isomorphism $\alpha$.

We have to show that  $\iota$ maps the line set  of
$\Aa(F_u,K)$ onto the set of all those intersections of $B$ with
lines of $\Aa(K,R)$  that contain at least two points.

Obviously, $\iota$ preserves collinearity.
Now consider three (different)
points  $k_0^\alpha a$, $k_1^\alpha a$, $k_2^\alpha a$
of $B$ that are collinear
in $\Aa(K,R)$. Then  $k_0^\alpha a=xk_1^\alpha a+(1-x)k_2^\alpha a$
holds for some
$x\in K$.
We compute $(k_0-k_2)^\alpha a=x(k_1-k_2)^\alpha a$. Since $a\in R^*$,
this means $x=((k_0-k_2)(k_1-k_2)^{-1})^\alpha \in K\cap K^\alpha=F_{u^{-1}}$.
Hence  $k_0,k_1,k_2$
are collinear in $\Aa(F_u,K)$.
\qed

A special case is the following:

\begin{exa}\label{Baer}
Let $R$ be a quaternion skew field and let $K$ be one of its maximal
commutative subfields. Then $\dim_K(R)=2$, i.e., $\Aa(K,R)$ is an  affine
plane. All lines of $\Aa(K,R)$ appear as blocks compatible at $\infty$ with
the standard block (since $R^*=R\setminus\{0\}$).
By \cite{havl-97}, Theorem 2,  the other elements of $\fC_\infty$ are
{\em affine Baer subplanes} of $\Aa(K,R)$, i.e., isomorphic to the
affine plane $\Aa(Z,K)$ over the center $Z$ of $R$, which in this case
coincides with the field $F_u$ for each $u\in R\setminus Z$.
\end{exa}


\begin{thebibliography}{10}

\bibitem{benz-73}
W.~Benz.
\newblock {\em Vorlesungen \"uber Geometrie der Algebren}.
\newblock Springer, Berlin, 1973.

\bibitem{blunck-99b}
A.~Blunck.
\newblock Reguli and chains over skew fields.
\newblock {\em Beitr\"age Algebra Geom.}, 41:7--21, 2000.

\bibitem{cohn-95}
P.M. Cohn.
\newblock {\em Skew Fields}.
\newblock Cambridge University Press, Cambridge, 1995.

\bibitem{havl-93b}
H.~Havlicek.
\newblock On the geometry of field extensions.
\newblock {\em Aequationes Math.}, 45:232--238, 1993.

\bibitem{havl-94a}
H.~Havlicek.
\newblock Spreads of right quadratic skew field extensions.
\newblock {\em Geom. Dedicata}, 49:239--251, 1994.

\bibitem{havl-97}
H.~Havlicek.
\newblock Affine circle geometry over quaternion skew fields.
\newblock {\em Discrete Math.}, 174:153--165, 1997.

\bibitem{herz-92}
A.~Herzer.
\newblock On sets of subspaces closed under reguli.
\newblock {\em Geom.~Dedicata}, 41:89--99, 1992.

\bibitem{herz-95}
A.~Herzer.
\newblock Chain geometries.
\newblock In F.~Buekenhout, editor, {\em Handbook of Incidence Geometry}.
  Elsevier, Amsterdam, 1995.

\bibitem{kroll-92}
H.-J. Kroll.
\newblock Zur {D}arstellung der {U}nterr\"aume von {K}ettengeometrien.
\newblock {\em Geom. Dedicata}, 43:59--64, 1992.

\bibitem{maeu+m+n-80}
H.~M\"aurer, R.~Metz, and W.~Nolte.
\newblock Die {A}utomorphismengruppe der {M}\"obiusgeometrie einer
  {K}\"orpererweiterung.
\newblock {\em Aequationes Math.}, 21:110--112, 1980.

\end{thebibliography}

\bigskip

Institut f\"ur Geometrie\\
der Technischen Universit\"at\\
Wiedner Hauptstra{\ss}e 8--10\\
A--1040 Wien\\
Austria
\end{document}